\documentclass[11pt]{article}

\usepackage[dvips]{graphicx}
\usepackage{amssymb}
\usepackage{amsmath}
\usepackage{amsthm}
\usepackage{latexsym}
\usepackage[english]{babel}
\usepackage{appendix}

\usepackage[usenames,dvipsnames]{color}

\newtheorem{theorem}{Theorem}[section]

\newtheorem{proposition}{Proposition}[section]

\def \C {\mathbb{C}}
\def \R {\mathbb{R}}

\def \P {\mathcal{P}}

\def\Xint#1{\mathchoice
{\XXint\displaystyle\textstyle{#1}}%
{\XXint\textstyle\scriptstyle{#1}}%
{\XXint\scriptstyle\scriptscriptstyle{#1}}%
{\XXint\scriptscriptstyle\scriptscriptstyle{#1}}%
\!\int}
\def\XXint#1#2#3{{\setbox0=\hbox{$#1{#2#3}{\int}$ }
\vcenter{\hbox{$#2#3$ }}\kern-.6\wd0}}

\def\dashint{\Xint-}

\begin{document}

\title{Functional Determinant on Pseudo-Einstein 3-manifolds}

\author{Ali Maalaoui$^{(1)}$}
\addtocounter{footnote}{1}
\footnotetext{Department of mathematics and natural sciences, American University of Ras Al Khaimah, PO Box 10021, Ras Al Khaimah, UAE. E-mail address:
{\tt{ali.maalaoui@aurak.ac.ae}}}

\date{}
\maketitle

\vspace{5mm}

{\noindent\bf Abstract} {Given a three dimensional pseudo-Einstein CR manifold $(M,T^{1,0}M,\theta)$, we establish an expression for the difference of determinants of the Paneitz type operators $A_{\theta}$, related to the problem of prescribing the $Q'$-curvature, under the conformal change $\theta\mapsto e^{w}\theta$ with $w\in \P$ the space of pluriharmonic functions. This generalizes the expression of the functional determinant in four dimensional Riemannian manifolds established in \cite{BO2}. We also provide an existence result of maximizers for the scaling invariant functional determinant as in \cite{CY}.}

\vspace{5mm}

\noindent

\vspace{4mm}

\noindent

\vspace{4mm}


\section{Introduction  and statement of the results}
There has been extensive work on the study of spectral invariants of differential operators defined on a Riemannian manifold $(M,g)$ and the relations to their conformal invariants, see for instance \cite{BO1,BO2,BO3} and the references therein. As an example, if we consider the two dimensional case with the pair of the Laplace operator $-\Delta_{g}$, and the associated invariant which is the scalar curvature $R_{g}$, we know that under conformal change of the metric $g\mapsto \tilde{g}=e^{2w}g$, one has the relation
$$R_{\tilde{g}}e^{2w}=-\Delta_{g}w+R_{g}.$$
It is also known that the spectrum of $-\Delta_{g}$ is discrete and can be written as $0< \lambda_{1}\leq \lambda_{2}\leq \cdots$ and the corresponding zeta function is then defined by $$\zeta_{-\Delta_{g}}(s)=\sum_{j=1}^{\infty}\frac{1}{\lambda_{j}^{s}}.$$
This series converges uniformly for $s>1$ and can be extended to a meromorphic function in $\C$ with $0$ as a regular value. The determinant of the operator $-\Delta_{g}$ can then be written as
$$det(-\Delta_{g})=e^{-\zeta'_{-\Delta_{g}}(0)}.$$
The celebrated Polyakov formula \cite{Pol1,Pol2}, states that if $\tilde{g}=e^{2w}g$  then
$$\ln\Big(\frac{det(-\Delta_{\tilde{g}})}{det(-\Delta_{g})}\Big)=\frac{-1}{12\pi}\int_{M}|\nabla w|^{2}+2Rw \ dv_{g},$$ 
for metrics with the same volume. The scaling invariant functional determinant $F_{2}$ can then be written as
$$F_{2}(w)=\frac{-1}{12\pi}\Big(\int_{M}|\nabla w|^{2}+2Rw \ dv_{g}-\Big(\int_{M}R\ dv_{g}\Big)\ln\Big(\dashint_{M} e^{2w}\ dv_{g}\Big)\Big).$$
Notice that the right hand side is a familiar quantity. It is the Beckner-Onorfi energy, \cite{Beck} and we know that
$$\int_{S^{2}}|\nabla w|^{2}+2Rw \ dv_{g} -\Big(\int_{M}R\ dv_{g}\Big)\ln\Big(\dashint_{S^{2}} e^{2w}\ dv_{g}\Big)\geq 0.$$
This notion of determinant was extended to dimension four for conformally invariant operators, keeping in mind that the substitute of the pair $(-\Delta_{g},R)$ in dimension two is the pair $(P_{g},Q_{g})$ in dimension four, where $P_{g}$ is the Paneitz operator and $Q_{g}$ is the Riemannian Q-curvature (\cite{BO2, Mal}). In addition, two new terms appear in the scaling invariant functional determinant expression. Indeed, if $(M,g)$ a $4$-dimensional manifold and $A_{g}$ a non-negative self-adjoint conformally covariant operator then, there exists $\beta_{1}, \beta_{2}$ and $\beta_{3} \in \R$ such that the scaling invariant functional determinant $F_{4}$ reads as
\begin{equation}\label{detr}
F_{4}(w)=\beta_{1}I(w)+\beta_{2}II(w)+\beta_{3}III(w),
\end{equation}
where 
$$\left\{\begin{array}{lll}
I(w)=4\int_{M}w|W_{g}|^{2}\ dv_{g}-(\int_{M}|W_{g}|^{2}\ dv_{g})\ln\Big(\dashint_{M}e^{4w}\ dv_{g}\Big),\\
\\
II(w)=\int_{M}wP_{g}w+4Q_{g}w\ dv_{g}-(\int_{M}Q_{g}\ dv_{g})\ln\Big(\dashint_{M}e^{4w}\ dv_{g}\Big),\\
\\
III(w)=12\int_{M}(\Delta_{g}w+|\nabla w|^{2})^{2}\ dv_{g}-4\int_{M}w\Delta_{g} R_{g}+R_{g}|\nabla w|^{2}\ dv_{g}.
\end{array}
\right.
$$
In the case of the sphere $S^{4},$ we see that the second term $II$ corresponds again to the four dimensional Beckner-Onofri energy. The existence and uniqueness of maximizers of this expression was heavily investigated and we refer the reader to \cite{CY, GM,Mal, Oki} and the references therein for an in depth study of this functional in the Riemannian setting.\\

Now let us move to the CR setting. We consider a 3-dimensional CR manifold $(M, T^{1,0}M,J,\theta)$ and we recall that the substitute for the pair $(P_{g},Q_{g})$ is $(P_{\theta}, Q_{\theta})$ where $P_{\theta}$ is the CR Paneitz operator and $Q_{\theta}$ is the CR Q-curvature, \cite{FH,GG}. The problem with this pair is that the total Q-curvature is always zero. In fact in pseudo-Einstein manifolds the Q-curvature vanishes identically. Hence, we do not have a decent substitute for the CR Beckner-Onofri inequality. Fortunately, if we restrict our study to pseudo-Einstein manifolds and variations in the space of pluriharmonic functions $\P$, then we have a better substitute for the pair $(P_{g},Q_{g})$ namely $(P'_{\theta},Q'_{\theta})$. These quantities were first introduced on odd dimensional spheres in \cite{Bran} and then on pseudo-Einstein manifolds in \cite{CaYa1,CaYa, Hir}. In particular one has a Beckner-Onofri type inequality involving the operator $P'_{\theta}$ acting on pluriharmonic functions as proved in \cite{Bran}. We also recall that the total $Q'$-curvature corresponds to a geometric invariant, namely the Burns-Epstein invariant $\mu(M)$ ( \cite{BE,CL}).\\
One is tempted to see what the spectral invariants of the operator $P'$ are or the restriction of $P'$ to the space $\P$ of pluri-harmonic functions and link them to geometric quantities such as the total Q'-curvature.\\
We recall that the quantity $Q'_{\theta}$ changes as follows: if $\tilde{\theta}=e^{w}\theta$ with $w\in \P$, then
\begin{equation}\label{qq}
P'_{\theta}w+Q'_{\theta}=Q'_{\tilde{\theta}}e^{2w}+\frac{1}{2}P_{\theta}(w^{2}),
\end{equation}
which we can write as
$$P'_{\theta}w+Q'_{\theta}=Q'_{\tilde{\theta}}e^{2w} \text{ mod } \P^{\perp}.$$
We let $\tau_{\theta}:L^{2}\to \P$  be the orthogonal projection on $\P$ with respect to the inner product induced by $\theta$ and set $A_{\theta}=\tau_{\theta} P'_{\theta}\tau_{\theta}$. Then equation $(\ref{qq})$ can be rewritten as
$$A_{\theta}w+\tau_{\theta}(Q'_{\theta})=\tau_{\theta}(Q'_{\tilde{\theta}}e^{2w}).$$
Prescribing the quantity $\overline{Q}'_{\theta}=\tau_{\theta}(Q'_{\theta})$ was thoroughly investigates in \cite{M1,CaYa1,QQ} mainly because of the property that
$$\int_{M}\overline{Q}'_{\theta}\ d\nu_{\theta}=\int_{M}Q'_{\theta}\ d\nu_{\theta}=-\frac{\mu(M)}{16\pi^{2}}.$$
We recall that in \cite{M}, we proved that the dual of the Beckner-Onofri inequality, namely the logarithmic Hardy-Littlewood-Sobolev inequality can be linked to the regularized zeta function of the operator $A_{\theta}$ evaluated at one. This was proved in the Riemannian setting in \cite{Mor1,Mor2, Ok1}.\\
In this paper, we will generalize the expression $(\ref{detr})$ by studying the determinant of the operator $A_{\theta}$. In all that follows we assume that $(M,T^{1,0}M,J,\theta)$ is an embeddable pseudo-Einstein manifold such that $P'_{\theta}$ is non-negative and has trivial kernel. First we show that

\begin{theorem}[Conformal Index]
Let $\zeta_{A_{\theta}}$ be the spectral zeta function of the operator $A_{\theta}$. Then $\zeta_{A_{\theta}}(0)$ is a conformal invariant in $\P$. Moreover,
$$\zeta_{A_{\theta}}(0)=\frac{-1}{24\pi^{2}}\int_{M}Q'_{\theta}\ dv_{\theta}-1.$$
\end{theorem}
In order to compute the determinant of the operator $A_{\theta}$ we introduce the quantities $\tilde{A}_{1}(w)$, $\tilde{A}_{2}(w)$ and $\tilde{A}_{3}(w)$ defined by
\begin{equation}
\left\{\begin{array}{lll}
\tilde{A}_{1}(w):=\int_{M}wA_{\theta}w+2Q'_{\theta}w\ d\nu-\frac{1}{c_{1}}\ln\Big(\dashint_{M}e^{2w}\ d\nu\Big),\\
\\
\tilde{A}_{2}(w):=2\int_{M}R\Big(\Delta_{b}w+\frac{1}{2}|\nabla_{b}w|^{2}\Big)-\Big(\Delta_{b}w+\frac{1}{2}|\nabla_{b} w|^{2}\Big)^{2}\ d\nu,\\
\\
\tilde{A}_{3}(w):=2\int_{M}w_{0}R-\frac{1}{3}w_{0}|\nabla_{b}w|^{2}-w_{0}\Delta_{b}w \ d\nu.
\end{array}
\right.
\end{equation}
One can also write $\tilde{A}_{2}(w)$ as
$$\tilde{A}_{2}(w)=2\int_{M}R\Big(\frac{\Delta_{b}e^{\frac{1}{2}w}}{2e^{\frac{1}{2}w}}\Big)-\Big(\frac{\Delta_{b}e^{\frac{1}{2}w}}{2e^{\frac{1}{2}w}}\Big)^{2}\ d\nu.$$
Then we have the following
\begin{theorem}
There exists $c_{2}$ and $c_{3}\in \R$ such that for all $w\in \P$, we have
\begin{equation}\label{det}
\ln\left(\frac{det(A_{\theta})}{det(A_{e^{w}\theta})}\right)=-\frac{1}{24\pi^{2}}\tilde{A}_{1}(w)+c_{2}\tilde{A}_{2}(w)-c_{3}\tilde{A}_{3}(w).
\end{equation}
\end{theorem}
Notice that the expressio $(\ref{det})$ is not scaling invariant, that is because for $\tilde{\theta}=c^{2}\theta$, with $c$ a positive constant, we have
 $$det(A_{\tilde{\theta}})=c^{-4\zeta_{\theta}(0)}det(A_{\theta}).$$
So we fix the volume $V$ of $(M,\theta)$ and define the scaling invariant functional
$$S_{A_{\theta}}=(\frac{Vol(\theta)}{V})^{\zeta_{A_{\theta}}(0)}det(A).$$
Now we can define the scaling invariant functional $F:W^{2,2}_{H}(M)\cap\P \to \R$ by
$$F(w)=\ln(S_{A_{\theta}})-\ln(S_{A_{e^{w}\theta}}),$$
where $W^{2,2}_{H}(M)$ is Folland-Stein space. Then one can write the following expression of $F$.
$$F(w)=c_{1}II(w)+c_{2}III(w)+c_{3}IV(w)$$
where
\begin{equation}\label{syst}
\left \{ \begin{array}{lllll}
II(w)=\int_{M}wA_{\theta}w+2Q'_{\theta}w \ d\nu -\int_{M}Q'_{\theta}\ d\nu \ln\Big(\dashint e^{2w} d \nu\Big).\\
\\
III(w)=\tilde{A}_{2}(w).\\
\\
IV(w)=-\tilde{A}_{3}(w).
\end{array}
\right.
\end{equation}
Notice that the functional $II$ is the CR Beckner-Onofi functional studied first in \cite{Bran} on the standard sphere $S^{3}$. In particular one has on the standard sphere
$$II(w)\geq 0.$$
 This functional was also investigated in \cite{CaYa1} and its critical points correspond t0 the pseudo-Einstein structures with constant $\overline{Q'}$-curvature. The functional $III$ is also similar to the Riemannian one defined in \cite{CY} and its critical points are pseudo-Einstein contact forms $\tilde{\theta}$ satisfying
$$\tilde{\tau}(\tilde{\Delta_{b}}\tilde{R})=0.$$
The functional $IV$ is a bit different, in fact if we let $\mathcal{H}$ defined by
$$\mathcal{H}(w)= R_{,0}-\frac{1}{3}|\nabla_{b}w|^{2}_{,0}-\frac{2}{3}div_{b}(w_{,0}\nabla_{b}w)+\Delta_{b}w_{,0}-(\Delta_{b}w)_{,0},$$
then the critical points of $IV$ satisfy, 
$$\tilde{\tau}(e^{-2w}\mathcal{H}(w))=0.$$
We set $a=\frac{\int_{M}Q'_{\theta}\ d\nu}{16\pi^{2}}$. Then, we show the following for the functional $F$:
\begin{theorem}
Assume that $c_{2}>0$ and $c_{3}\geq 0$. Then there exists a constant $\mu$ depending on $\theta$, such that if
\begin{equation}\label{cond}
c_{3}<\mu\Big(\sqrt{25c_{2}^{2}+\frac{1}{3\pi^{2}}c_{2}(1-a)}-5c_{2}\Big),
\end{equation}
then $F$ has a maximizer $w_{\infty}\in W^{2,2}_{H}(M)\cap \P$ under the constraint $\int_{M}e^{2w}dv=1$. Moreover, this maximizer satisfies the Euler-Lagrange equation
$$\tau_{\theta_{\infty}}\Big[-\frac{1}{24\pi^{2}}\tilde{Q'_{\theta}}+c_{2}\tilde{\Delta_{b}}\tilde{R}+c_{3}e^{-2w}\mathcal{H}(w)\Big]=cte,$$
where the tilde here, refers to quantities computed using the contact form $\theta_{\infty}=e^{w_{\infty}}\theta$.
\end{theorem}
Notice the condition $(\ref{cond})$ implies in particular that $\int_{M}Q'_{\theta} d\nu<16\pi^{2}$. Hence, as a consequence, we have that $(M,T^{1,0}M,\theta)$ is not equivalent to the standard sphere. We point out that in order to verify the sharpness of $(\ref{cond})$ one needs to check specific examples which is not as easy as in the Riemannian case, since we are dealing with CR pluriharmonic functions and we do lack explicit examples of manifolds where one can have an explicit expression of the spectrum of the $P'$-operators.\\

\noindent
{\bf Acknowledgement}
The author wants to acknowledge that the problem treated in this paper was initiated by a question of Carlo Morpurgo. The author is grateful for the fruitful discussions that led to the finalized version of this paper. \\

\section{Heat Coefficients and Conformal Invariance}
Let $(M,T^{1,0}M,\theta)$ be a pseudo-Einstein 3-manifold and $P'_{\theta}$ its $P'$-operator defined by
\begin{equation}\label{Pan}
P_\theta^\prime f = 4\Delta_b^2 f - 8 \text{Im}\left(\nabla^1(A_{11}\nabla^1 f)\right) - 4\text{Re}\left(\nabla^1(R\nabla_1 f)\right).
\end{equation}
Denote by $\tau:L^{2}(M)\to \P$ the orthogonal projection on the space of pluriharmonic functions with respect to the $L^{2}$-inner product induced by $\theta$. We consider the operator $A_{\theta}=\tau P'_{\theta} \tau$ and for the conformal change $\tilde{\theta}=e^{w}\theta$, with $w\in \P$, we let $$A_{\tilde{\theta}}=\tau_{\tilde{\theta}}\Big(e^{-2w}A_{\theta}\Big),$$
where $\tau_{\tilde{\theta}}$ is the orthogonal projection with respect to the $L^{2}$-inner product induced by $\tilde{\theta}$.\\
In order to evaluate and manipulate the spectral invariants, we need to study the expression of the heat kernel of the operator $A_{\theta}$. Unfortunately, this operator is not elliptic or sub-elliptic (as an operator on $C^{\infty}(M)$), and does not have an invertible principal symbol in the sense of $\Psi_{H}(M)$-calculus (see \cite{P}). In fact $A_{\theta}$ can be seen as a Toeplitz operator, and one might adopt the approach introduced in $\cite{Bout}$ in order to study it. But instead, we will modify the operator in order to be able to use the classical computations done for the heat kernel.\\
Consider the operator $\mathcal{L}=A_{\theta}+\tau^{\perp}L\tau^{\perp}$, where $L$ is chosen so that $\mathcal{L}$ has an invertible principal symbol in $\Psi_{H}^{4}(M)$. Notice that $\tau \mathcal{L}=\mathcal{L}\tau=A_{\theta}$. Based on \cite{P}, if $\mathcal{K}$ is the heat kernel of $\mathcal{L}$ one has the following expansion near zero:
$$\mathcal{K}(t,x,x)\sim \sum_{j=0}^{\infty}\tilde{a}_{j}(x)t^{\frac{j-4}{4}}+\ln(t)\sum_{j=1}^{\infty}t^{j}\tilde{b}_{j}(x).$$
Since $e^{-t\mathcal{L}}=e^{-tA_{\theta}}\tau +e^{-tL}\tau^{\perp}$, we have that the kernel $K$ of $e^{-tA_{\theta}}$ which is the restriction to $\P$ of $\mathcal{K}$, reads as
\begin{equation}\label{exp}
K(t,x,x)\sim\sum_{j=0}^{\infty}t^{\frac{j-4}{4}}a_{j}(x)+\ln(t)\sum_{j=1}^{\infty}t^{j}b_{j}(x),
\end{equation}
and this will be the main expansion that we will be using for the rest of the paper.\\
Now we want to define the infinitesimal variation of a quantity under a conformal change. Fix $w\in \P$ and for a given quantity $F_{\theta}$ depending on $\theta$ denote $\delta F_{\theta}:=\frac{d}{dr}_{|r=0}F_{e^{rw}\theta}.$ 
Next, we define the zeta function of $A_{\theta}$ by
$$\zeta_{A_{\theta}}(s):=\sum_{j=1}^{\infty}\frac{1}{\lambda_{j}^{s}}.$$
where $0<\lambda_{1}\leq \lambda_{2} \leq \cdots$, is the spectrum of the operator $A_{\theta}:W^{2,2}(M)\cap \P \to \P$. In what follows $TR[A]$ is to be understood as the trace of the operator $A$ in $\P$. Then we have the following proposition.
\begin{proposition}
With the notations above, we have
$$\zeta_{A_{\theta}}(0)=\int_{M}a_{4}(x)dx-1.$$
Moreover,
$$\delta \zeta_{A_{\theta}}(0)=0$$
and
$$\delta \zeta_{A_{\theta}}'(0)=2\int_{M}w(a_{4}(x)-\frac{1}{V})\ d\nu,$$
where $V=\int_{M}d\nu_{\theta}$ is the volume of $M$.
\end{proposition}
{\it Proof:}
Most of the computations in this part are relatively standard and they can be found in \cite{BO1,BO2,BO3} in the Riemannian setting. First we use the Mellin transform  and $(\ref{exp})$ to write
\begin{align}
\zeta_{A_{\theta}}(s)&=\frac{1}{\Gamma(s)}\int_{0}^{\infty}t^{s-1}(TR[e^{-tA_{\theta}}]-1)\ dt\notag\\
&=\frac{1}{\Gamma(s)}\Big(-\frac{1}{s}+\int_{0}^{1}t^{s-1}\sum_{j=0}^{N}t^{\frac{j-4}{4}}\int_{M}a_{j}(x)\ d\nu \ dt +\int_{M}t^{s-1}O(t^{\frac{N+1-4}{4}})\ dt\notag\\
&+\sum_{j=1}^{N}\int_{0}^{1}t^{j+s-1}\ln(t)\int_{M}b_{j}\ d\nu\ dt+ \int_{0}^{1}t^{s-1}O(t^{N+1}\ln(t))\ dt +\int_{1}^{\infty}t^{s-1}\sum_{j=1}^{\infty}e^{-\lambda_{j} t}\ dt\Big)\notag\\
&=\frac{1}{\Gamma(s)}\Big(\frac{-1}{s}+\sum_{j=0}^{N}\frac{1}{s+\frac{j-4}{4}}\int_{M}a_{j}(x)\ d\nu +\int_{0}^{1}t^{s-1}O(t^{\frac{N+1-4}{4}})\ dt\notag\\
&+\sum_{j=1}^{N}\frac{1}{(s+j)^{2}}\int_{M}b_{j}\ d\nu+\int_{0}^{1}t^{s-1}O(t^{N+1}\ln(t))\ dt\int_{1}^{\infty}t^{s-1}\sum_{j=1}^{\infty}e^{-\lambda_{j} t}\ dt\Big).\notag
\end{align}
Since, $\Gamma$ has a simple pole at $s=0$ with residue 1, we see that by taking $s\to 0$, there are only two terms that survive, leading to
$$\zeta_{A_{\theta}}(0)=\int_{M}a_{4}(x)\ d\nu-1.$$
Next we move to the study of the variation of $\zeta_{A_{\theta}}$. Let $f\in C^{\infty}(M)$ and $v\in \P$. Then
$$\int_{M}\tau_{rw}(f)v\ d\nu_{rw}=\int_{M}fve^{2rw}\ d\nu.$$
Differentiating with respect to $r$ and evaluating at $0$ yields
$$\int_{M}(\delta\tau(f)+2w\tau(f)-2wF)v\ d\nu=0.$$
Hence,
$$\delta \tau (f)=2\tau(wf-w\tau(f)).$$
If we let $M_{w}$ be the multiplication by $w$ then
$$\delta \tau =2(\tau M_{w}-\tau M_{w}\tau).$$
In particular, if $f\in \P$ then $\delta \tau (f)=0$. \\
Next we want to evaluate $\delta A_{\theta}$. Recall that $A_{e^{rw}\theta}=\tau_{e^{rw}\theta}e^{-2rw}A_{\theta}$. Therefore,
$$\delta A_{\theta} = \delta\tau A_{\theta}-2\tau M_{w}A_{\theta}=-2\tau M_{w} A_{\theta}.$$
Thus,
\begin{align}
\delta TR[e^{-tA_{\theta}}]&=-tTR[\delta A_{\theta}e^{-tA_{\theta}}]\notag\\
&=2tTR[\tau M_{w} A_{\theta}e^{-tA_{\theta}}]=2tTR[M_{w}A_{\theta}e^{-tA_{\theta}}]\notag.
\end{align}
The last equality follows from $TR[AB]=TR[BA]$ and $e^{-tA_{\theta}}\tau =e^{-tA_{\theta}}$.
But
$$-tM_{w}A_{\theta}e^{-tA_{\theta}}=\frac{d}{d\varepsilon}_{|\varepsilon=0}M_{w}e^{-t(1+\varepsilon)A_{\theta}}$$
Using the expansion $(\ref{exp})$, we have
$$K(t(1+\varepsilon),x,x)\sim\sum_{j=0}^{\infty}(1+\varepsilon)^{\frac{j-4}{4}}t^{\frac{j-4}{4}}a_{j}(x)+H,$$
where $H$ is the logarithmic part. Hence, comparing the terms in the expansion after integration, we get:
$$\delta \int_{M}a_{j}\ d\nu=\frac{4-j}{2}\int_{M}wa_{j}\ d\nu.$$
In particular, we have $\delta\int_{M}a_{4}\ d\nu=0$.\\
Similarly
$$
\Gamma(s)\zeta_{A_\theta}(s)=\Gamma(s)(\zeta_{A_{\theta}}(0)+s\zeta'_{A_{\theta}}(0)+O(s^{2})).
$$
Hence, since $\delta \zeta_{A_{\theta}}(0)=0$, and $s\Gamma(s)\sim 0$ when $s\to 0$, we  have
$$\delta \zeta'_{A_{\theta}}(0)=[\Gamma(s)\delta \zeta_{A_{\theta}}(s)]_{s=0}.$$
But,
\begin{align}
\Gamma(s)\delta \zeta_{A_{\theta}}(s)&=\int_{0}^{\infty}2t^{s}TR[wA_{\theta}e^{-tA_{\theta}}]\ dt\notag\\
&=-\int_{0}^{\infty}2t^{s}\frac{d}{dt}TR[we^{-tA_{\theta}}]\ dt\notag\\
&=\int_{0}^{\infty}2st^{s-1}TR[w(e^{-tA_{\theta}}-\frac{1}{V})]\ dt.
\end{align}
Using again the expansion $(\ref{exp})$ and a similar computation as in the previous case, yields
$$\delta \zeta_{A_{\theta}}'(0)=2\int_{M}w(a_{4}-\frac{1}{V})\ d\nu.$$
\hfill$\Box$

\begin{proposition}
There exists $c\not=0$ such that
$$\zeta_{A_{\theta}}(0)=c\int_{M}Q'_{\theta}\ d\nu-1.$$
Moreover $c=-\frac{1}{24\pi^{2}}$.
\end{proposition}
{\it Proof:}
 First we notice that $a_{4}$ is a pseudo-Hermitian invariant of order $-2$, that is $$a_{4,e^{r}\theta}=e^{-2r}a_{4,\theta},$$
for all $r\in \R$. So from \cite{Hir}, we have the existence of $c_{1},c_{2},c_{3},c_{4},c_{5}\in \R$ such that
$$a_{4}=c_{1}Q'_{\theta}+c_{2}\Delta_{b}R+c_{3}R_{,0}+c_{4}R^{2}+c_{5}Q_{\theta},$$
where $Q'_{\theta}= 2\Delta_b R - 4 |A|^2 + R^2$ and $Q_{\theta}=-\frac{2}{3}\Delta_{b}R+2Im( A_{11,\bar{1}\bar{1}})$.
Since we are in a pseudo-Einstein manifold and $w\in \P$ we can assume that $Q_{\theta}=0$. So after integration, we have
$$\int_{M}a_{4}\ d\nu=c_{1}\int_{M}Q'_{\theta}\ d\nu+c_{4}\int_{M}R^{2}\ d\nu.$$
Since $\int_{M}a_{4}\ d\nu$ is invariant under the conformal change $e^{w}\theta$, it is easy to see that $c_{4}=0$. Hence,
$$\int_{M}a_{4}\ d\nu=c_{1}\int_{M}Q'_{\theta}\ d\nu.$$
Next we want to calculate $c_{1}$ (compare to \cite{ST}, where the invariant $k_{2}$ is always 0). We take the case of the sphere $S^{3}$. Based on the computations in \cite{Bran}, we have
$$\zeta_{A_{\theta}}(s)=2\sum_{j=1}^{\infty}\frac{j+1}{(j(j+1))^{s}}=2\sum_{j=2}^{\infty}\frac{1}{j^{2s-1}}\Big(\frac{1}{1-\frac{1}{j}}\Big)^{s}.$$
Using the expansion of $(1-\frac{1}{j})^{-s}=1+\frac{s}{j}+\frac{s(s+1)}{2j^2}+sO(\frac{1}{j^{3}})$, we see that
$$\zeta_{A_{\theta}}(s)=2(\zeta_{R}(2s-1)-1+s(\zeta_{R}(2s)-1)-\frac{s(s+1)}{2}(\zeta_{R}(2s+1)-1))+sH(s).$$
with $H(s)$ holomorphic near $s=0$ and $\zeta_{R}$ the classical Riemann Zeta function. Now we recall that $\zeta_{R}$ is regular at $s=0$ and $s=-1$ but has a simple pole at $s=1$ with residue equal to $1$. Hence
$$\zeta_{A_{\theta}}(0)=2(-\frac{1}{12}-1+\frac{1}{4})=-\frac{5}{3}\not=0.$$
Knowing that $\int_{S^{3}}Q'_{\theta}\ d\nu=16\pi^{2}$ we have 
$$16\pi^{2}c-1=-\frac{5}{3}.$$
Thus,
$$c=-\frac{1}{24\pi^{2}}.$$
\hfill$\Box$

\section{The expression for the Determinant}
Recall that in the previous section, we found that $a_{4}=c_{1}Q'+c_{2}\Delta_{b}R+c_{3}R_{,0}$. In particular
\begin{align}
\delta \zeta_{A_{\theta}}'(0)&=\int_{M}2w\Big(a_{4}(x)-\frac{1}{V}\Big)\ d\nu\notag\\
&=c_{1}\int_{M}2w\Big(Q'_{\theta}-\frac{1}{c_{1}V}\Big)\ d\nu+c_{2}\int_{M}2R\Delta_{b}w\ d\nu -c_{3}\int_{M}2w_{,0}R\ d\nu\notag\\
&=c_{1}A_{1}+c_{2}A_{2}+c_{3}A_{3}.\notag
\end{align}
We will calculate the change of each term under conformal change of $\theta$. The easiest term to handle is the first one.
Indeed, recall that if $\tilde{\theta}=e^{w}\theta$ then
$$\tilde{Q'_{\theta}}e^{2w}=P'_{\theta}w+Q'_{\theta} \text{ mod } \P^{\perp},$$
$$\tilde{R}=\Big[R-|\nabla_{b}w|^{2}-2\Delta_{b}w\Big]e^{-w},$$
and
$$\tilde{\Delta}_{b}f=e^{-w}\Big[\Delta_{b}f+\nabla_{b}f\cdot\nabla_{b}w\Big].$$
So if $\hat{\theta}_{u}=e^{uw}\theta$, we have
$$\int_{M}2w\Big[\hat{Q'_{\theta}}-\frac{1}{c_{1}}\frac{1}{\hat{V}}\Big]\ d\hat{v}=\int_{M}2uwP'_{\theta}w+2Q'w-\frac{1}{c_{1}}\frac{2we^{2uw}}{\int_{M}e^{2uw}\ d\nu}\ d\nu.$$
Integrating $u$ in $[0,1]$ yields
$$\tilde{A}_{1}(w)=\int_{M}wA_{\theta}w+Q'_{\theta}w\ d\nu-\frac{1}{c_{1}}\ln(\dashint_{M}e^{2w}\ d\nu).$$
For the second term, we have
\begin{align}
\int_{M}\hat{R}\hat{\Delta}_{b}w\ d\hat{v}&=\int_{M}\Big[R-u^{2}|\nabla_{b}w|^{2}-2u\Delta_{b}w\Big]\Big[\Delta_{b}w+u|\nabla_{b}w|^{2}\Big]\ d\nu\notag\\
&=\int_{M}R\Delta_{b}w-u^{2}|\nabla_{b}w|^{2}\Delta_{b}w-2u(\Delta_{b}w)^{2}+Ru|\nabla_{b}w|^{2}-u^{3}|\nabla_{b}w|^{4}\notag\\
&\quad-2u^{2}|\nabla_{b}w|^{2}\Delta_{b}w\ d\nu. \notag
\end{align}
In particular after integrating over $u$ between $0$ and $1$, we get
\begin{align}
\tilde{A}_{2}(w)&=2\int_{M}R\Delta_{b}w-|\nabla_{b}w|^{2}\Delta_{b}w-(\Delta_{b}w)^{2}+\frac{1}{2}R|\nabla_{b}w|^{2}-\frac{1}{4}|\nabla_{b}w|^{4}\ d\nu\notag\\
&=2\int_{M}R\Delta_{b}w-\Big(\Delta_{b}w+\frac{1}{2}|\nabla_{b} w|^{2}\Big)^{2}+\frac{R}{2}|\nabla_{b}w|^{2}\ d\nu. \notag
\end{align}
Next we compute
$$\int_{M}\hat{T}w \hat{R}d\hat{v}=\int_{M}\Big[w_{,0}R-u^{2}w_{,0}|\nabla_{b}w|^{2}-2uw_{,0}\Delta_{b}w\Big]\ d\nu,$$
where $T$ is the characteristic vectorfield of $\theta$ and we are adopting the notatoion $Tf=f_{,0}$. Integrating as above yields:
$$\tilde{A}_{3}(w)=2\int_{M}w_{,0}R-\frac{1}{3}w_{,0}|\nabla_{b}w|^{2}-w_{.0}\Delta_{b}w\ d\nu.$$
Therefore, one has
$$\zeta'_{\tilde{A_{\theta}}}(0)-\zeta'_{A_{\theta}}(0)=c_{1}\tilde{A}_{1}(w)+c_{2}\tilde{A}_{2}(w)-c_{3}\tilde{A}_{3}(w)$$
or equivalently
$$\ln\Big(\frac{\det(A_{\theta})}{\det(A_{\tilde{\theta}})}\Big)=c_{1}\tilde{A}_{1}(w)+c_{2}\tilde{A}_{2}(w)-c_{3}\tilde{A}_{3}(w).$$
\hfill$\Box$

\section{Scaling invariant functional}
We focus now on the study of the functional $F$ defined by
$$F(w)=c_{1}II(w)+c_{2}III(w)+c_{3}IV(w),$$
where $c_{1}=-\frac{1}{24\pi^{2}}$. For the sake of notation, we will keep using $c_{1}$ instead of its numerical value. We will also be using constants $C_{k}$ depending on $M$ and $\theta$.\\
We recall first that there exists a constant $C$ such that
\begin{equation}\label{bo}
\frac{1}{16\pi^{2}}\int_{M}wA_{\theta}w +2Q'w\ d\nu -\ln\Big(\dashint e^{2w} \ d\nu \Big)\geq C.
\end{equation}
In fact this follows form the CR Beckner-Onofri inequality proved in \cite{Bran} and also treated in\cite{CaYa}.
Since the functional $F$ is scaling invariant (that is $F(w+c)=F(w)$, we can assume without loss of generality that $\overline{w}=\int_{M}w \ d\nu=0$. So we can write $F$ as
\begin{align}
F(w)&=2\int_{M}a_{4}w\ d\nu+c_{1}\Big(\int_{M}wA_{\theta}w \ d\nu-\int_{M}Q'_{\theta}\ d\nu\ln\Big(\dashint e^{2w}\Big)\Big)\notag\\
&\quad-2c_{2}\int_{M}\Big(\Delta_{b}w+\frac{1}{2}|\nabla_{b}w|^{2}\Big)^{2}\ d\nu +c_{2}\int_{M}R|\nabla_{b}w|^{2}\ d\nu \notag\\
&\quad+c_{3}\int_{M} w_{,0}\Big(\frac{1}{3}|\nabla w|^{2}+\Delta_{b}w\Big)\ d\nu.\notag
\end{align}
Using $(\ref{bo})$, we have
$$\int_{M}Q'_{\theta}\ d\nu \ln\Big(\dashint e^{2w}\ d\nu\Big)\leq \frac{\int_{M}Q'_{\theta}\ d\nu}{16\pi^{2}}\Big[\int_{M}wA_{\theta}w+2Q'_{\theta}w \ d\nu\Big]-C.$$ 
Therefore, for $a=\frac{\int_{M}Q'_{\theta}\ d\nu}{16\pi^{2}}$, we have
\begin{align}
c_{1}\Big(\int_{M}wA_{\theta}w \ d\nu-\int_{M}Q'_{\theta}\ d\nu\ln\Big(\dashint e^{2w}\Big)\Big) &\leq c_{1}\Big[(1-a)\int_{M}wA_{\theta}w\ d\nu +2a\int Q'_{\theta}w \ d\nu\Big] +C_{1}\notag\\
&\leq 4c_{1}(1-a)\int_{M}(\Delta_{b}w)^{2} \ d\nu +C_{2}\int_{M}|\nabla_{b}w|^{2} \ d\nu\notag\\
&\quad +2ac_{1}\int Q'_{\theta}w \ d\nu +C_{3},\label{eq0}
\end{align}
where in the second line we used the expression $(\ref{Pan})$. On the other hand, for the mixed term of $III(w)$, we have for every $\alpha>0$,
\begin{equation}\label{eq1}
2\int_{M}\Delta_{b}w|\nabla_{b}w|^{2}\ d\nu \leq \alpha \int_{M} (\Delta_{b}w)^{2}+\frac{1}{\alpha}\int_{M}|\nabla_{b}w|^{4}\ d\nu.\end{equation}
Next, we let $\lambda$ denote the best constant appearing in the estimate
$$\|f_{,0}\|_{L^{2}}\leq \lambda \|\Delta_{b} f\|_{L^{2}}.$$
Then we have
\begin{align}
2\int_{M}w_{,0}\Big(\frac{1}{3} |\nabla_{b}|^{2}+\Delta_{b}w\Big) \ d\nu &\leq 2\Big(\lambda \|\Delta_{b}w\|_{L^{2}}\|\frac{1}{3}|
\nabla_{b}w|^{2}+\Delta_{b}w\|_{L^{2}}\Big)\notag\\
&\leq 2\lambda\Big(\|\Delta w\|_{L^{2}}^{2}+\frac{1}{3}\|\Delta w\|_{L^{2}}\||\nabla_{b}w|^{2}\|_{L^{2}}\Big)\notag\\
&\leq \lambda \Big((2+\frac{\alpha}{3})\|\Delta w\|_{L^{2}}^{2}+\frac{1}{3\alpha}\int_{M}|\nabla_{b}w|^{4}\ d\nu\Big).\label{eq2}
\end{align}
Hence, combining $(\ref{eq0})$, $(\ref{eq1})$ and $(\ref{eq2})$ and assuming that $c_{2}>0$ and $c_{3}\geq 0$, we get
\begin{align}
c_{1}II(w)+c_{2}III(w)+c_{3}IV(w)\leq& \Big(4c_{1}(1-a)+c_{2}(\alpha-2)+c_{3}(2+\frac{\alpha}{3})\lambda\Big)\int_{M}(\Delta_{b}w)^{2}\ d\nu \notag\\
&+\Big(c_{2}(\frac{1}{\alpha}-\frac{1}{2})+c_{3}\frac{\lambda}{3\alpha}\Big)\int_{M}|\nabla_{b}w|^{4}\ d\nu +C_{4}\int_{M}|\nabla_{b}w|^{2}\ d\nu\notag\\
&+2\int_{M}a_{4}w\ d\nu+2ac_{1}\int_{M}Q'_{\theta}w\ d\nu+C_{5}.\label{eq4}
\end{align}

Now we need to choose $\alpha$ in a way that the coefficients of $\int_{M}(\Delta_{b}w)^{2}\ d\nu$ and $\int_{M}|\nabla_{b}w|^{2}\ d\nu$ are both negative. For this to happen, one needs that

$$\left\{\begin{array}{lll}
4c_{1}(1-a)-2c_{2}+2\lambda c_{3}<0,\\
\alpha<\frac{2c_{2}-2\lambda c_{3}-4c_{1}(1-a)}{c_{2}+\frac{1}{3}\lambda c_{3}},\\
\frac{1}{\alpha}<\frac{c_{2}}{2c_{2}+\frac{2}{3}\lambda c_{3}}.
\end{array}
\right.$$
Hence, we need
$$\frac{2c_{2}+\frac{2}{3}\lambda c_{3}}{c_{2}}<\frac{2c_{2}-2\lambda c_{3}-4c_{1}(1-a)}{c_{2}+\frac{1}{3}\lambda c_{3}}.$$
This is possible if condition $(\ref{cond})$ is satisfied for $\mu=\frac{\lambda}{3}$. This yields in particular that if $w_{k}\in W^{2,2}_{H}(M)\cap \P$ is a maximizing sequence for $F$, then there exists $C>0$ such that
$$\int_{M}(\Delta_{b}w_{k})^{2}\ d\nu+\int_{M}|\nabla_{b}w_{k}|^{4} \ d\nu\leq C.$$
Indeed, we have $$-\varepsilon<F(0)\leq F(w_{k}).$$
Therefore, there exists $c>0$ such that
\begin{align}
-\varepsilon \leq& -c\Big(\int_{M}(\Delta_{b}w_{k})^{2}\ d\nu+\int_{M}|\nabla_{b}w_{k}|^{4}\ d\nu\Big)+C_{4}\int_{M}|\nabla_{b}w_{k}|^{2}\ d\nu\notag\\
&+2\int_{M}a_{4}w_{k}\ d\nu +2ac_{1}\int_{M}Q'_{\theta}w_{k}\ d\nu +C_{5}.\notag
\end{align}
Thus
$$\int_{M}(\Delta_{b}w_{k})^{2}\ d\nu+\int_{M}|\nabla_{b}w_{k}|^{4}\ d\nu\leq C_{6}\Big(\int_{M}|\nabla_{b}w_{k}|^{4}\ d\nu \Big)^{\frac{1}{2}}+C_{7},$$
where here we used the fact that $\overline{w}=0$ and
$$\int_{M}fw \ d\nu\leq \|f\|_{L^{2}}\|w\|_{L^{2}} \leq \|f\|_{L^{2}}\|\nabla_{b}w\|_{L^{2}}\leq \|f\|_{L^{2}}V^{\frac{1}{4}}\|\nabla_{b}w\|_{L^{4}}.$$
Hence, $(w_{k})_{k}$ is bounded in $W^{2,2}_{H}(M)$ and have a weakly convergent subsequence, that converges to $w_{\infty}$. Passing to the $\limsup$ in $F(w_{k})$, we get that the weak limit $w_{\infty}$ is in fact a maximizer of $F$.\\
Finally, based on the remark below $(\ref{syst})$, we see that the critical points of $F$ under the constraint $\int_{M}e^{2w}dv=1$ satisfy the equation
$$\tilde{\tau}\Big[c_{1}\tilde{Q'_{\theta}}+c_{2}\tilde{\Delta_{b}}\tilde{R}+c_{3}e^{-2w}\mathcal{H}\Big]=cte.$$

\end{document}